\documentstyle[12pt]{article}

\begin{document}

\thispagestyle{empty}
\setcounter{page}{0}

\vspace*{-8ex}
 \null \hfill MPI-PhT/95-1 \\
 \null \hfill GK-MP 9502/17 \\
 \null \hfill q-alg/9502007 \\
 \null \hfill February 1995 \\[15mm]

\begin{center}

{\bf \LARGE
  Dynamical Symmetries in q-deformed\\  [1ex]
  Quantum Mechanics} \\  [6ex]

  A. Lorek$^1$ and J. Wess$^{1,2}$
   \\ [5ex]

 $^1$ Max-Planck-Institut f\"{u}r Physik \\
   Werner-Heisenberg-Institut\\
 F\"{o}hringer Ring 6 , D - 80805 M\"{u}nchen, Germany \\
 Tel. (89) 32308-1, Fax (89) 3226704\\ [3ex]

 $^2$ Sektion Physik, Universit\"{a}t M\"{u}nchen \\
   Theresienstr.\,37, D - 80333 M\"{u}nchen, Germany \\ [20ex]

\end{center}

\begin{abstract}

\noindent
The dynamical algebra of the q-deformed harmonic oscillator is constructed.
As a result, we find the free deformed Hamiltonian as well as the Hamiltonian
of the deformed oscillator as a complicated, momentum dependent interaction
Hamiltonian in terms of the usual canonical variables.\\
Furthermore we construct a well-defined algebra SU$_q$(1,1) with consistent
conjugation properties and comultiplication. We obtain non lowest weight
representations of this algebra.

\end{abstract}

\clearpage

\renewcommand{\theequation}{\arabic{section}.\arabic{equation}}
\newcommand{\qi}{\frac{1}{q}}
\newcommand{\qiq}{\frac{1}{q^2}}
\newcommand{\nn}{\nonumber}
\newcommand{\be}{\begin{equation}}
\newcommand{\ee}{\end{equation}}
\newcommand{\bea}{\begin{eqnarray}}
\newcommand{\eea}{\end{eqnarray}}
\newcommand{\ZZ}{\mbox{Z\hspace{-1.35mm}Z}}
\newcommand{\NN}{\mbox{I\hspace{-0.58mm}N}}
\newcommand{\RR}{\mbox{I\hspace{-0.62mm}R}}
\newcommand{\CC}{\mbox{C\hspace{-2.0mm}l}}
\newcommand{\Unity}{\mbox{1\hspace{-1.1mm}I}}

\section{Introduction}

\setcounter{equation}{0}

Dynamical groups provide a useful tool for solving dynamical systems.
For the harmonic oscillator it is the SU(1,1) group that has the Hamiltonian
as one generator while the other generators generate the spectrum.
Diagonalizing the Hamiltonian can be treated as a problem of group
representations. \\
It is natural to follow a similar strategy for the q-deformed
harmonic oscillator. We base this oscillator and consequently the dynamical
group on q-deformed canonical variables. The generators are quadratic
expressions in terms of these variables, it is easy to see that they form a
q-deformed algebra which is a deformation of the U(1,1) algebra. It is
possible to find a central piece in this q-deformed algebra as well and to
reduce the algebra to a three generator algebra. It is, however, very
difficult to formulate the conjugation properties that would justify
to call it a SU$_q$(1,1) algebra. This is much easier in the four generator
version - especially as we can start from the conjugation properties of the
canonical variables.\\
The thus derived conjugation properties of the algebra are generalized
to the case of all real values of the Casimir
operator whereas the dynamical group algebra
is characterized by a particular value of the Casimir.\\[4mm]
We obtain a general definition of a U$_q$(1,1) algebra
\bea
B \,A - q^4 A\, B &=& - i \,(1+q^2) \,u\, C \nn\\
C \,A - q^2 A\, C &=& - i q \,(1+q^2) (1+q^4) \,u\, A \nn \\
C \,B - \qiq B\, C &=& \frac{i}{q^3} \,(1+q^2) (1+q^4) \,u\, B \\[1mm]
u\, A &=& \qiq \,A\, u\nn \\
u\, B &=& q^2\,B\, u \nn\\
u\, C &=& C\, u\nn
\eea
\\[1mm]
Comultiplication is an algebraic property, it can be abstracted from the
SU$_q$(2) version of the algebra. One possibility, consistent with the
algebra (1.1) is
\bea
\Delta(A) & = & A \otimes u + u\,\tau_3^{\frac{1}{2}}\otimes A \nn \\
\Delta(B) & = & B \otimes u + u\,\tau_3^{\frac{1}{2}}\otimes B  \\
\Delta(C) & = & C \otimes u + u\,\tau_3 \otimes C
+\frac{q(q^8-1)}{\gamma}\, \tau_3^{\frac{1}{2}} \, B \otimes A
+\frac{q(q^8-1)}{q^4\gamma}\,\tau_3^{\frac{1}{2}} \, A \otimes B \nn \\
\Delta(u) & = &  u \otimes u \nn\\[1mm]
\tau_3 & = & u^{-2}\,V^2\,\left(
1\,-\,\frac{(1-q^4)(1-q^8)}{2\,q^2}\cal{C}\right) \nn\\
\cal{C} &=& \left[ \frac{2}{q^2 \alpha\beta} A\, B + \frac{2}{q^4 \gamma
(q^2+\frac{1}{q^2})} C\; (\frac{1}{\gamma} C + q^2 u) \right] V^{-2} \nn\\
V &=&  u - \frac{1}{\gamma} (q^2 - \frac{1}{q^2})\; C \nn\\
\gamma &=& \frac{i}{q}\,(1+q^2)(1+q^4) \nn
\eea
\\[2mm]
The conjugation properties are:
\bea
\overline{A} & = & A \hspace{20mm}\overline{B} = B \nn \\
\overline{C} & = & \left[ 1\, - \, \frac{(1-q^4)(1-q^8)}{2\,q^2}\, \cal{C}
\right]^{-\frac{1}{2}}\left[ C\,-\,\frac{i}{2\,q}\, (1-q^8)(1+q^2)(1+q^4)
\,V\,\cal{C}\right] \nn\\
\overline{u} & = & \left[ 1\, - \, \frac{(1-q^4)(1-q^8)}{2\,q^2}\, \cal{C}
\right]^{\frac{1}{2}} V \, = \, \tau_3^{\frac{1}{2}}\,u
\eea
\\[2mm]
Representations of this algebra where A or B are diagonal are constructed
for values of the Casimir ${\cal C}  <  \frac{2\,q^2}{(1-q^4)(1-q^8)}$.
These are not
lowest weight representations.\\
Representations of the algebra where $A+B$ are diagonal which correspond
to the harmonic oscillator are not easy to obtain, we retreat to perturbation
theory.\\

Finally we use the fact that the generators of the q-deformed algebra are in
the enveloping algebra of the undeformed one. This allows us to express the
Hamiltonian of the q-deformed system in terms of the common Heisenberg
variables $x$, $p_x$ with $[p_x,x]= -i$. It is a Hamiltonian with a
complicated, momentum dependent interaction:
\bea
{\cal H} \:=\: \frac{1}{2}\,B &=& 2\;\; \mbox{\large{\it p$_x$}} \;\;
\sqrt{\frac{q+q^{-1} - 2 \cos((x p_x+p_x x)\,h)}{(q^2-q^{-2})^2 \,
((x p_x+p_x x)^2+1)}}\;\;\: \mbox{\large{\it p$_x$}} \nn\\
q &=& e^h
\eea
\\
This is the free Hamiltonian of a q-deformed system with the spectrum
\be
E_n = \frac{1}{2} \,\pi_0^2 \,q^{2n}\;\;, \;\;\;\;\;\;\;n\in \NN
\ee
The system has no lowest energy state. $\pi_0$ is a constant characterizing
the representation.\\

This serves to demonstrate that systems with complicated interactions as (1.4)
can have a q-deformed algebra that governs their kinematics. This kinematics
can then be used to find solutions of the interacting system. The interactions
turn out to be momentum dependent, q-deformation gives us a handle to study
the properties of such interactions more systematically.\\[2mm]
It is also interesting that this kinematics in a natural way leads to space
and momentum variables that have discrete spectra - space appears to
be quantized.

\section{Dynamical Variables and Algebras}

\setcounter{equation}{0}

The dynamical group algebras will be based on the q-deformed canonical
variables as they were introduced in ref.\cite{schw}:
\bea
p\, \xi - q\, \xi \,p & = & - i q \,u \nn\\
u \,p & = & q\, p\, u \\
u \,\xi & = & \qi\, \xi\, u \nn
\eea
These variables are subject to conjugation properties, consistent with
the algebra:
\be
\overline{p} = p \;\;,\hspace{1cm}\overline{\xi} = \xi \;\;,\hspace{1cm}
\overline{u} = \qi \,u^{-1}
\ee
Certain polynomials in the variables p and $\xi$ (together with u) form
algebras,
they are of the SU$_q$(1,1) type.\\
The first example are the homogeneous polynomials of second degree:
\be
A = \xi^2 \;\;,\hspace{1cm} B = p^2 \;\;,\hspace{1cm} C = p \xi + q^3 \xi p
\ee
which form the algebra:
\bea
B \,A - q^4 A\, B &=& - i \,(1+q^2) \,u\, C \nn\\
C \,A - q^2 A\, C &=& - i q \,(1+q^2) (1+q^4) \,u\, A \nn \\
C \,B - \qiq B\, C &=& \frac{i}{q^3} \,(1+q^2) (1+q^4) \,u\, B \\[1mm]
u\, A &=& \qiq \,A\, u\nn \\
u\, B &=& q^2\,B\, u \nn\\
u\, C &=& C\, u\nn
\eea
Another example in analogy to the SU(1,1) subalgebra of the Virasoro
algebra is:\\
\be
L_0 = \xi p \;\;,\hspace{1cm}L_+ = p\;\;,\hspace{1cm}L_- = \xi^2 p
\ee
\\
with the corresponding algebra
\bea
L_0 L_+ - \qi L_+ L_0&=&i u L_+ \nn \\
L_- L_0 - \qi L_0 L_-&=&i q u L_-  \\
L_+ L_- - q^2 L_- L_+&=&-i q (1+q^2) u L_0 \nn \\[2mm]
u\, L_0 &=& L_0\, u \nn\\
u\, L_+ &=& q \,L_+\, u \nn\\
u\, L_- &=& \qi\, L_-\, u \nn
\eea
An algebra of the same type can be constructed in terms of q-deformed
variables in three dimensions \cite{weich}:
\bea
A' &=& \vec{P}^2 \nn \\
B' &=& \vec{X}^2 \\
C' &=& \Lambda^{-\frac{1}{2}} \frac{1}{q^2-1} \left(
(1+q^2)^2 W^2 - q^4 (q^8+1) \Lambda - 2 q^2 \right) \nn
\eea
Using formulas of ref.\cite{weich} we obtain:
\bea
B' A' - q^8 A' B' &=& - \frac{(q^2+1)(q^4+1)}{4\, q^6}\Lambda^{\frac{1}{2}} C'
\nn\\
C' A' - q^4 A' C' &=& q^4 (q^8+1)(q^4+1)(q^2+1) \Lambda^{\frac{1}{2}} A' \\
C' B' - \frac{1}{q^4} B' C' &=& -\frac{1}{q^4} (q^8+1)(q^4+1)(q^2+1)
\Lambda^{\frac{1}{2}} B' \nn \\[2mm]
\Lambda\; A' &=& \frac{1}{q^8} A'\; \Lambda \nn\\
\Lambda\; B' &=& q^8 B'\; \Lambda \nn\\[1mm]
\Lambda\; C' &=& C'\; \Lambda \nn
\eea
\\
As it is well known a SU$_q$(1,1) algebra can also be realized in terms of
deformed harmonic oscillator variables \cite{mcf}
\be
a a^+ - \frac{1}{r} a^+ a = 1
\ee
Here we define
\bea
W_- &=& \hat{\alpha}^{-1} a^2 \nn\\
W_+ &=& \hat{\beta}^{-1} (a^+)^2 \\
W_0 &=& \hat{\gamma}^{-1} \left( 1 + \frac{1}{r} (1+\frac{1}{r})a^+ a \right)
\nn
\eea
and choose the normalization
\be
\hat{\gamma} = (1+r)(1+\frac{1}{r^2})\;\;,\hspace{1cm}
\hat{\alpha} \hat{\beta} = - r (1+r) \hat{\gamma}
\ee
to obtain the algebra
\bea
r\; W_0 W_+ - \frac{1}{r}\; W_+ W_0 &=& W_+ \nn\\
r\; W_- W_0 - \frac{1}{r}\; W_0 W_- &=& W_- \\
\frac{1}{r^2}\; W_+ W_- - r^2\; W_- W_+ &=& W_0 \;\;\;\;\; .\nn
\eea
This algebra (which we will refer to as W-algebra) has been studied
by Curtright and Zachos \cite{curt}. By
relating it to the SU$_q$(2) algebra they succeeded in finding a Casimir
operator and opened the way to derive comultiplication rules.

Our strategy is to relate to this W-algebra the algebras (2.4), (2.6) and (2.8)
and to obtain this way Casimir and comultiplication rules for all these
algebras.\\[2mm]
The respective relations are:\\[2mm]
Algebra (2.4):\\[-3mm]
\bea
W_+ = \frac{1}{\beta} \,u^{-1}\, B & \hspace{2cm}&\gamma = \frac{i}{q}
(1+q^2)(1+q^4) \nn\\
W_- = \frac{1}{\alpha} \,u^{-1}\, A& &\alpha \,\beta = \frac{1}{q^3}
(1+q^2)^2(1+q^4) \\
W_0 = \frac{1}{\gamma} \,u^{-1}\, C & & r = q^2 \nn
\eea
As usual, the algebra leaves a scaling freedom for $W_+$ and $W_-$, therefore
only the product $\alpha \beta$ can be fixed.\\[3mm]
Algebra (2.6):\\[-3mm]
\bea
W_+ = \frac{1}{\beta'} \,u^{-1}\, L_+ &\hspace{2cm}&\gamma' = i\; q \nn\\
W_- = \frac{1}{\alpha'} \,u^{-1}\, L_- & &\alpha' \,\beta' = q (1+q^2) \\
W_0 = \frac{1}{\gamma'} \,u^{-1}\, L_0 & & r = q \nn
\eea
\\
Algebra (2.8):\\[-3mm]
\bea
W_+ = \frac{1}{\alpha''} \,\Lambda^{-\frac{1}{2}}\, A'& \hspace{2cm} &\gamma''
= (q^8+1)(q^4+1)(q^2+1) \nn\\
W_- = \frac{1}{\beta''} \,\Lambda^{-\frac{1}{2}} \, B'&  &\alpha''
\,\beta'' = \frac{q^8-1}{4 (q^2-1) q^{10}} \gamma'' \\
W_0 = \frac{1}{\gamma''} \,\Lambda^{-\frac{1}{2}} \, C'& & r =
\frac{1}{q^4} \nn
\eea
\\
These identifications are unique up to the symmetry of the algebra (2.12)
under the transformation:\\[-8mm]
\bea
W_0 &\rightarrow& - W_0 \nn\\
W_+ &\rightarrow& W_- \\
r &\rightarrow& \frac{1}{r} \nn
\eea
\\
The expression of the Casimir for the W-algebra (2.12) is\\[-2mm]
\bea
\cal{C} & = &\left[ 2 W_-\; W_+ + \frac{2}{r^2 (r+\frac{1}{r})} W_0
(W_0+r) \right] \; \left[1 - (r - \frac{1}{r}) W_0 \right]^{-2}\;\;\; .
\eea
With the identifications given above it can be used to find the Casimir
for all these algebras. As a special example we give the expression for
the algebra (2.4):
\bea
\cal{C} &=& \left[ \frac{2}{q^2 \alpha\beta} A\, B + \frac{2}{q^4 \gamma
(q^2+\frac{1}{q^2})} C\; (\frac{1}{\gamma} C + q^2 u) \right] V^{-2}\\
V &=&  u - \frac{1}{\gamma} (q^2 - \frac{1}{q^2})\; C \nn
\eea
\\[3mm]
If we now express the elements of the algebra in terms of the dynamical
variables we obtain:\\[2mm]
Algebra (2.4)\\[-6mm]
\bea
\cal{C} & = & -\, \frac{2(1-q^6)}{(1-q^2)(1+q^4)(1+q^2)^4}
\eea
Algebra (2.6)\\[-6mm]
\bea
\cal{C} & = & 0 \nn
\eea
Algebra (2.8)\\[-2mm]
\bea
\cal{C} & = & - \frac{1}{2}\,\frac{q^{12}}{1+q^8} + \frac{q^{12}(q^8+1)\;W^2}
{((q^2+1)^2 W^2-2 q^2)^2 (q^4+1)^2} \left[ q^4 (q^2+1)^2 L^2 +1 \right]
\eea
Observe that the value of the Casimir in the three dimensional case
depends on the angular momentum which is already true in the classical
case. W is related to the angular momentum via the relation:
\be
W^2 - 1 = q^4 (q^2-1)^2\; L^2
\ee
For the algebra (2.12) in terms of the harmonic oscillator variables we obtain
the same value for $\cal{C}$ as in (2.19) replacing $q^2$ by $r$.\\

\section{Comultiplication and conjugation properties}

\setcounter{equation}{0}

As was mentioned above Curtright and Zachos \cite{curt} have studied the
algebra (2.12) and have related it to the usual SU$_q$(2) algebra. If the
parameter r of (2.12) and the parameter q of SU$_q$(2) are identified
($r=q$), the relation is particularly simple
\bea
W_+ & = & \sqrt{\frac{q}{q+{\displaystyle \qi}}} \: T_+ \nn \\
W_- & = & \qi \,\sqrt{\frac{q}{q+{\displaystyle \qi}}} \: T_- \\
W_0 & = & \frac{q}{q^2-{\displaystyle \qiq}}\,(1-\tau_3)
\,-\,\qi\,\frac{q-{\displaystyle \qi}}{q+{\displaystyle \qi}}\:T_+\,T_- \nn
\eea
A straightforward calculation starting from (2.12) easily reproduces the T
algebra ($\lambda:= q - {\displaystyle \qi}$):
\bea
\qi\,T_+\,T_-\,-\,q\,T_-\,T_+ & = & T_3 \nn \\
q^2 \,T_3\,T_+\,-\,\qiq\,T_+\,T_3 & = & (q+\qi)\,T_+ \nn \\
q^2 \,T_-\,T_3\,-\,\qiq\,T_3\,T_- & = & (q+\qi)\,T_- \\[1mm]
\tau_3 & := & 1-\lambda T_3 \nn\\
\tau_3\,T_+ & = &\frac{1}{q^4}\,T_+\,\tau_3  \nn\\
\tau_3\,T_- & = & q^4\,T_-\,\tau_3\nn
\eea
\\
As the comultiplication rule for this algebra is well known it is easy
to give such a rule for the W-algebra:
\bea
\Delta(W_{\pm}) & = & W_{\pm} \otimes 1 + \tau_3^{\frac{1}{2}} \otimes
W_{\pm} \nn \\
\Delta(W_0) & = & W_0 \otimes 1 + \tau_3 \otimes W_0 -
q\,\lambda\,\tau_3^{\frac{1}{2}}\,W_+ \otimes W_-
- \qi\,\lambda\,\tau_3^{\frac{1}{2}}\,W_- \otimes W_+ \\
 \tau_3 & = & 1 - \frac{\lambda}{q}\,(q+\qi)\,W_0 -
\frac{\lambda^2}{q^2}\,(q+\qi)\,W_+\,W_- \nn
\eea
It can be directly verified that this comultiplication rule is
compatible with the W-algebra (2.12).

If we define u to have a grouplike comultiplication rule with $\alpha$
being a free parameter
\be
\Delta(u)  =  \alpha \: \left( u \otimes u \right)
\ee
then using (2.13) we obtain a consistent comultiplication scheme for the
algebra (2.4):
\bea
\Delta(A) & = & \alpha \left( \,A \otimes u +
u\,\tau_3^{\frac{1}{2}}\otimes A \right) \nn \\
\Delta(B) & = & \alpha \left( \,B \otimes u +
u\,\tau_3^{\frac{1}{2}}\otimes B \right) \\
\Delta(C) & = & \alpha \left( \,C \otimes u + u\,\tau_3 \otimes C
+\frac{q(q^8-1)}{\gamma}\, \tau_3^{\frac{1}{2}} \, B \otimes A
+\frac{q(q^8-1)}{q^4\gamma}\,\tau_3^{\frac{1}{2}} \, A \otimes B \right) \nn \\
\tau_3 & = & u^{-2}\,V^2\,\left(
1\,-\,\frac{(1-q^4)(1-q^8)}{2\,q^2}\cal{C}\right) \nn
\eea
This $\tau_3$ is the same as in (3.2) and (3.3), expressed in the
variables of the algebra (2.4) with V and $\cal{C}$ defined in
(2.18).\\
To verify this comultiplication rule for the algebra (2.4) is more tedious.
To guess it seems to be quite a task.\\

To construct representations for the algebras and for physical
applications conjugations properties are essential. One way to find
conjugation rules in our case is to look at the algebra as it is expressed in
dynamical variables.\\
For the algebra (2.4) this leads to:
\bea
\overline{A} & = & A \hspace{15mm}\overline{B} = B \hspace{2cm} \overline{q} =
q
\nn \\
\overline{C} & = & \qi\, \frac{1+q^4}{1+q^2}\,C
\,+\,i\,\frac{1-q^6}{1+q^2}\,u \nn \\
\overline{u} & = & \qi\, \frac{1+q^4}{1+q^2}\,u
\,-\,\frac{i}{q^2}\,\frac{1-q^2}{1+q^2}\,C  \\
\cal{C} & = & -\,\frac{2\,(1-q^6)}{(1-q^2)(1+q^4)(1+q^2)^4}
\,=\,\overline{\cal C} \nn
\eea
This conjugation rule, however, is only compatible with the algebra
(2.4) for the particular value of the Casimir given above.\\
A general conjugation rule will have to involve the Casimir operator
itself.
To find such a conjugation rule we start from $\overline{A}=A$ and
$\overline{B}=B$. $\overline{C}$ and $\overline{u}$ can then be
derived by conjugating the algebra relations (2.4). The result is:
\bea
\overline{A} & = & A \hspace{20mm}\overline{B} = B \nn \\
\overline{C} & = & \left[ 1\, - \, \frac{(1-q^4)(1-q^8)}{2\,q^2}\, \cal{C}
\right]^{-\frac{1}{2}}\left[ C\,-\,\frac{i}{2\,q}\, (1-q^8)(1+q^2)(1+q^4)
\,V\,\cal{C}\right] \nn\\
\overline{u} & = & \left[ 1\, - \, \frac{(1-q^4)(1-q^8)}{2\,q^2}\, \cal{C}
\right]^{\frac{1}{2}} V \, = \, \tau_3^{\frac{1}{2}}\,u
\eea
The elements V and $\tau_3$ were already defined in (2.18) and
(3.5) respectively.
The Casimir turns out to be hermitean:
\be
\overline{\cal{C}} = \cal{C}
\ee

For $q=1$ the algebra with its conjugation properties reduces to SU(1,1)
together with a central hermitean element $u$. We can, as suggested from
(2.2), further impose for $q=1$: $\overline{u} = u^{-1}$. This yields $u=\pm
1$. \\
The algebra (2.4) has for arbitrary $q$ a central element as well. This
element is $uV$, it can be verified that it commutes with all the elements
of the algebra. Under conjugation we find:
\be
u\,V\,=\,\overline{u\,V}\,=\,\overline{V}\,\overline{u}
\ee
\\
{}From the representations of the algebra we shall learn that a relation
\be
u^{-1}\,=\,\beta^2\,\overline{u}
\ee
($\beta$ a real parameter) can be imposed consistently. In this case the
central element becomes
\be
 u\,V\,= \frac{q^2-1}{\beta^2} \left[1\,-\,\frac{(1-q^4)(1-q^8)}{2\,q^2}\,
 {\cal C} \right]^{-\frac{1}{2}}
\ee

We now study the coalgebra under conjugation. For the algebra (2.4) with
the conjugation rule (3.7) we can write the comultiplication (3.4) and (3.5)
in a simpler form:
\bea
\Delta(A) & = & \alpha \left( A \otimes u + \overline{u} \otimes A \right)
\nn \\
\Delta(B) & = & \alpha \left( B \otimes u + \overline{u} \otimes B \right) \\
\Delta(C) & = & \alpha \left( C \otimes u + \frac{\overline{u}^2}{u} \otimes
C -i\,q^2\,(q^2-1)\,\frac{\overline{u}}{u} \,B \otimes A
 -i\,\qiq\,(q^2-1)\,\frac{\overline{u}}{u} \,A \otimes B \right) \nn \\
\Delta(u) & = & \alpha \left( u \otimes u \right) \nn
\eea
In this form it is easy to see that:
\bea
\overline{\Delta(A)} & = & \sigma \cdot \Delta_{\sigma}(\overline{A}) \nn \\
\overline{\Delta(B)} & = & \sigma \cdot \Delta_{\sigma}(\overline{B})
\eea
where $\sigma$ twists the factors in the direct product
\be
\sigma\, \cdot\, X\otimes Y = Y \otimes X
\ee
and $\Delta_{\sigma}$ is the twisted comultiplication
\bea
\Delta(M) &=& M_{(1)} \otimes M_{(2)} \nn\\
\Delta_{\sigma}(M) &=& M_{(2)} \otimes M_{(1)}
\eea
where the indices 1 and 2 may denote different representations for instance.
\\[2mm]
A short calculation shows that the same is true for $\Delta(C)$ and
$\Delta(u)$:
\bea
\overline{\Delta(C)} & = & \sigma \cdot \Delta_{\sigma}(\overline{C}) \nn \\
\overline{\Delta(u)} & = & \sigma \cdot \Delta_{\sigma}(\overline{u})
\eea
where $\Delta_{\sigma}(\overline{C})$ and $\Delta_{\sigma}(\overline{u})$
are defined by (3.7) and the fact that $\Delta_{\sigma}$
is an algebra homomorphism. \\
If in addition (3.10) is imposed to eliminate the central piece we find
$\alpha=\beta$ in order that
\be
\Delta(u^{-1}) = \beta^2\,\Delta(\overline{u})
\ee
\\
Thus we have found a q-deformation of SU(1,1) as a coalgebra as well.
This is not quite the Hopf algebra structure we want as it contains the twist
operator $\sigma$.\\

The algebra with its comultiplication and conjugation properties has an
interesting symmetry. The following transformation leaves the algebra
invariant:
\bea
q & \rightarrow & \qi \nn \\[1mm]
A & \rightarrow & \tilde{A} \,=\, \rho \, B \nn \\
B & \rightarrow & \tilde{B} \,=\, \eta \, A  \\
C & \rightarrow & \tilde{C} \,=\, \kappa C \nn \\
u & \rightarrow & \tilde{u} \,=\,- q^4 \kappa \, u
\hspace{15mm}\rho\,\eta = q^6 \kappa^2\nn
\eea
For real values of $\rho$, $\eta$ and $\kappa$ the conjugation properties
remain unchanged. For the comultiplication the parameter $\alpha$ in (3.12)
has to be scaled:
\be
\alpha \rightarrow \tilde{\alpha} \,=\,- \frac{1}{q^4\, \kappa}\, \alpha
\ee
If we eliminate the central piece as above by using (3.10) we obtain
\be
1=q^8\,\kappa^2 \tilde{\beta}^2 / \beta^2
\ee
This is consistent with (3.19) and $\alpha = \beta$.\\[3mm]
The Casimir as defined in (2.18) changes
\bea
\cal{C} & \rightarrow & q^8\,\cal{C}
\eea
However, the expression $q^4\cal{C}$ would not change and would serve as a
Casimir as well.

The symmetry (3.18) turns out to be a generalisation of the symmetry (2.16)
of the W-algebra, it includes the conjugation and comultiplication
as well.\\

\section{Representations}

\setcounter{equation}{0}

The representations of the algebra (2.4) with conjugation properties (3.7) and
B diagonal are easy to construct.\\
We first observe that a scaling of B$\rightarrow$tB and A$\rightarrow$t$^{-1}$A
does neither change the algebra (2.4) nor the Casimir (2.18). With this in mind
we write
\be
B\,|\nu>\,=\,b_{\nu} \,|\nu>
\ee
and assume that the eigenvalues are not degenerate. From the
commutation relation with u follows
\be
B\,u\,|\nu>\,=\,\qiq\,b_{\nu} \,u\,|\nu> \;\;\;.
\ee
With an eigenvalue of B all the $q^{2n}$ multiples are eigenvalues as well.
Therefore
\bea
B\,|\nu>\,&=&\,b_{0}\,q^{2\nu} \,|\nu> \nn \\
u\,|\nu>\,&=&\,\alpha_{\nu} \,|\nu -1>
\eea
As mentioned above $b_0$ cannot be determined from the algebra. We will see
that these states for a given $b_0$ are actually sufficient to construct a
representation of the algebra.

{}From the B,C relation follows
\be
(b_{\nu}-\qiq b_{\mu})\,<\mu|C|\nu> \,=\, \frac{i}{q^3} (1+q^2)(1+q^4)\;b_{\nu}
\,\alpha_{\nu}\, \delta_{\mu,\nu-1}
\ee
This shows that $<\mu|C|\nu>$ is zero except for $\mu=\nu-1$ or $\mu=\nu+1$.\\
The equation (4.4) determines the matix element for $\mu=\nu-1$:
\be
<\mu|C|\mu+1> = i\,q\,\frac{q^4+1}{q^2-1} \; \alpha_{\mu+1}
\ee
{}From the u,C relation follows:
\be
\alpha_{\mu+1}\,<\mu+1|C|\mu> = \alpha_{\mu} \,<\mu|C|\mu-1>
\ee
To get more information about the still undetermined matrix element of $C$
we have
to use the conjugation property of u in equation (3.7). The matrix element
of V is easily computed from the definition of V in (2.18):
\be
<\mu|V|\nu> = \frac{1-q^4}{q^2\,\gamma} \, \delta_{\mu,\nu+1}\,<\nu+1|C|\nu>
\ee
Now it is straightforward to calculate the matrix element of C using equation
(4.6) and assuming $\alpha_{\nu}$ to be real.\\
The result is
\be
<\mu|C|\nu> \,=\, i\,q\,\frac{q^4+1}{q^2-1} \, \alpha_{\nu}\,\delta_{\mu,\nu-1}
        - i\,q\,\frac{q^4+1}{q^2-1} \, \alpha_{\nu}\,
         \frac{\delta_{\mu,\nu+1}}{\sqrt{1-\frac{(1-q^4)(1-q^8)}{2\,q^2}
          \,\cal{C}}}
\ee
In this formula $\cal{C}$ stands for the eigenvalue of the Casimir.\\[2mm]
{}From the conjugation properties of C follows after some calculation and the
assumption that the square root in (4.8) is real
\be
\alpha_{\mu+1}\,=\,\alpha_{\mu}
\ee
For an imaginary square root there is no consistent solution.\\
Now the matrix elements of $C$ are determined up to one real parameter
$\alpha_0$.

We follow the same analysis for the operator A. From the A,u relation follows
\be
<\mu+1|A|\nu>\,=\,\qiq\,<\mu|A|\nu-1>
\ee
More information on these matrix elements is derived from the B,A relation:
\be
(b_{\mu}- q^4 b_{\nu})\,<\mu|A|\nu> = - i \,(1+q^2)\,
\alpha_{\mu+1}\,<\mu+1|C|\nu>
\ee
As the matrix elements of C are known we obtain for $\mu \neq \nu+2$ the
non vanishing matrix elements of A:
\bea
<\mu|A|\mu> &=& \frac{q}{(q^2-1)^2 b_{\mu}} \;
                \frac{1+q^4}{\sqrt{1-\frac{(1-q^4)(1-q^8)}{2\,q^2} \,\cal{C}}}
                \,\alpha_0^2 \nn \\
<\mu|A|\mu+2> &=& -\frac{q}{(q^2-1)^2 b_{\mu}} \;\alpha_0^2
\eea
The only other non vanishing matrix element of A is obtained from the
hermiticity of A:
\be
<\mu+2|A|\mu> = -\frac{q}{(q^2-1)^2 b_{\mu}} \,\alpha_0^2
\ee
The $\mu$-dependence via $b_{\mu}$ is consistent with (4.10).\\
The remaining A,C relation is identically statisfied.

If we finally impose the relation (3.10), eliminating the central piece of
the algebra, we can fix the constant $\alpha_0$ to be:
\be
\alpha_0 = \frac{1}{\beta}
\ee
This way we have obtained representations of the algebra for real
eigenvalues of the Casimir operator satisfying
\be
{\cal C}  <  \frac{2\,q^2}{(1-q^4)(1-q^8)}
\ee
which guarantees that the square root in (4.8) is real.\\[2mm]
For the realization of the algebra in terms of dynamical variables the
square root equals $\frac{1+q^4}{q\,(1+q^2)}$ and the eigenvalue of the
Casimir being negative lies well in the above range. In this case
we obtain a representation which follows directly from the
representation of the dynamical variables as it was constructed in ref.
\cite{schw} and \cite{heb}.\\

\section{SU$_q$(1,1) and the q-deformed harmonic oscillator}

\setcounter{equation}{0}

The group SU(1,1) is the dynamical group of the harmonic oscillator.\\
The generator $A+B$ is the Hamiltonian of the q-deformed case
as can be seen from (2.3). Representing the algebra SU$_q$(1,1) with $A+B$
diagonal would give information on the energy eigenvalues.\\
Unfortunately the construction of the representations with $A+B$ diagonal
is not as straightforward as the construction of the representations
in the previous chapter. Therefore we retreat to perturbation theory in q,
having the advantage that we can start from a discrete spectrum.

We first rewrite the algebra in terms of the generators X, Y:
\bea
X & = & A + B \nn \\
Y & = & A - B
\eea
With this definition (2.4) takes the following form
($\gamma = \frac{i}{q} (1+q^2)(1+q^4)$):
\bea
X\,Y - Y\,X + \frac{1-q^4}{1+q^4} (X^2-Y^2) & = &
-4i\,\frac{1+q^2}{1+q^4}\,u\,C
\nn \\
C\,X - q^2 X\,C + C\,Y - q^2 Y\,C & = & -q^2\,\gamma\,u\, (X+Y) \nn\\
C\,X - \qiq X\,C - (C\,Y - \qiq Y\,C) & = & \qiq\,\gamma\,u\, (X-Y) \nn
\eea
\bea
u\,X &=& \frac{1}{2} (q^2+\qiq)\,X\,u - \frac{1}{2} (q^2-\qiq)\,Y\,u \\
u\,Y &=& -\frac{1}{2} (q^2-\qiq)\,X\,u + \frac{1}{2} (q^2+\qiq)\,Y\,u \nn\\
u\,C &=& C\,u \nn
\eea
\\
In this version of the algebra for $q \ne 1$ the relations do not separate
into X and Y relations. This causes difficulties in constructing the
representation.\\[2mm]
The Casimir takes a particularly simple form in the X, Y, C generators:
\be
{\cal C}  =  \qi \left[ -X^2 + Y^2 + \frac{2}{q(1+q^4)}\,C^2 \right] \,
              \frac{1}{\gamma^2 V^2}
\ee
For $q=1$ the relations (5.2) and (5.3) reduce to
\bea
X_0 Y_0 - Y_0 X_0 &=& -4i\,C_0 \nn\\
C_0 X_0 - X_0 C_0 &=& -4i\,Y_0 \nn\\
C_0 Y_0 - Y_0 C_0 &=& -4i\,X_0 \\
\cal{C} &=& \frac{1}{16} (X_0^2 - Y_0^2 - C_0^2) \nn
\eea
These relations have the well known realization in terms of the creation and
annihilation operators of the undeformed harmonic oscillator.
\bea
X_0 &=& 2\,a^+ a + 1 \nn\\
Y_0 &=& -(a^2 + (a^+)^2) \\
C_0 &=& i\,(a^2-(a^+)^2)\nn\\
{\cal C}_0 &=& -\frac{3}{16} \nn
\eea
This suggests to introduce a linear combination of $Y$ and $C$ that reduces
to $a^2$ or $(a^+)^2$ for $q=1$:
\be
K^{\pm} = -\frac{1}{2} (Y \pm i C )
\ee
\\
In the basis of the harmonic oscillator ($|n> = \frac{1}{n!}(a^+)^n|0>$)
the matrix elements of $K_0^{\pm}$
simply are
\be
<n|K_0^+|n+2> \,=\, <n+2|K_0^-|n> \,=\, \sqrt{(n+1)(n+2)}
\ee
We expand around $q=1$ to first order
\bea
q &=& 1 + h \nn \\
u &=& 1 + h\, F \nn\\
X &=& X_0 + h\, X_1 \nn\\
Y &=& Y_0 + h\, Y_1 \\
C &=& C_0 + h\, C_1 \nn\\
K^{\pm} &=& K_0^{\pm} + h\, K_1^{\pm} \nn
\eea
{}From (2.2) and (3.6) we determine $F$ in terms of $C_0$ and learn about the
conjugation properties of $C_1$
\bea
F &=& - \frac{1}{2} (iC_0 +1) \;,\;\;\;\; \overline{F} = -F-1 \nn\\
\overline{C}_1 &=& C_1 - 3i
\eea
This yields the conjugation properties of $K_1^{\pm}$:
\be
\overline{K_1^+} = K_1^- + \frac{3}{2}
\ee
We expand the algebra (5.2) and take appropriate linear combinations to obtain
in first order in $h$ the following independent equations (after
substituting for F and using the Casimir in lowest order):
\bea
K_0^{\pm}X_1 - X_1 K_0^{\pm} + K_1^{\pm}X_0 - X_0 K_1^{\pm} &=&
   \pm 6 K_0^{\mp} \pm 4 K_1^{\pm} -3 \nn\\[2mm]
K_0^- K_1^+ - K_1^+ K_0^- + K_1^- K_0^+ - K_0^+ K_1^- &=& -3X_0 - 2X_1
\eea
The matrix elements of the first equation are:
\bea
&&\sqrt{(n+1)(n+2)}\,<n+2|X_1|m> - \sqrt{m(m-1)}\,<n|X_1|m-2> \nn\\
&=& (4+2(n-m))\, <n|K_1^+|m> - 3\,\delta_{n,m} + 6\sqrt{n(n-1)}\:\delta_{n,m+2}
\eea
For the special case $m=n+2$ we obtain
\be
<m|X_1|m> \,=\, <m-2|X_1|m-2>\, =:\, {\cal X}
\ee
The diagonal matrix elements of $X_1$ are constant.\\[2mm]
For $m \ne n+2$ we can solve equation (5.12) for the matrix elements of
$K_1^+$:
\bea
<n|K_1^+|m> &=& \frac{1}{4+2(n-m)} ( \sqrt{(n+1)(n+2)} \,<n+2|X_1|m>\\
&& -\sqrt{m(m-1)}\,<n|X_1|m-2> + 3\,\delta_{n,m} - 6\sqrt{n(n-1)}\,
\delta_{n,m+2} )\nn
\eea
\\
The corresponding matrix elements of $K_1^-$ are obtained analogously, this
is consistent with the conjugation (5.10).\\[3mm]
The second of the equations (5.11) has the following matrix elements
\bea
&&-\sqrt{n(n-1)}\,<n-2|K_1^+|m> + \sqrt{(m+2)(m+1)}\,<n|K_1^+|m+2>\nn\\
&&-\sqrt{m(m-1)}\,<n|K_1^-|m-2> + \sqrt{(n+1)(n+2)}\,<n+2|K_1^-|m> \nn \\[1mm]
&&= 3(2n+1)\,\delta_{n,m} + 2\,<n|X_1|m>
\eea
In the special case of $n=m$
we obtain a recursion formula for the matrix elements of
\mbox{$\sqrt{n(n-1)} \,2 \,Re <n-2|K_1^+|n>$}.
This recursion formula can be solved and we obtain:
\be
\sqrt{n(n-1)} \,2\, Re <n-2|K_1^+|n> = D_1 - \frac{1}{2}(3-2{\cal X})n
+ \frac{3}{2}n^2
\ee
$D_1$ is a free constant. It can be determined along with ${\cal X}$ by
looking at the equation for $n=0$ and $n=1$:
\be
D_1 = 0\; ,\hspace{7mm}{\cal X} = 0
\ee
If we assume consistently with (5.13) and (5.17) $X_1=0$ and in
addition \linebreak
Im$<n-2|K_1^+|n> = 0$, it is easy to obtain a
solution of all the equations to first order:
\bea
X &=& X_0 \nn\\
Y &=& Y_0 \nn\\
C &=& C_0 + \frac{3}{2}(i+C_0)\,h\\
u &=& 1 + \frac{i}{2}(i-C_0)\,h\nn
\eea
It is easy to verify directly that this is a solution starting from the
defining equations (5.4).

We finally study the behaviour of this solution under the reflection
(3.18). It leads to another solution
if we take the tilded variables as solution
of the algebra with $h$ replaced by $-h$ and compute the untilded
variables in terms of the tilded ones. The result is
\bea
X' &=& X_0 \nn\\
Y' &=& -Y_0 \nn\\
C' &=& -C_0 + \frac{3}{2}(i-C_0)\,h \\
u' &=& 1 + \frac{i}{2}(i+C_0)\,h \nn
\eea
This new solution can be obtained by a unitary transformation from the old
one, this can be
seen from the representation in terms of the harmonic oscillator operators
(5.5).\\

\section{q-Deformation and interaction}

\setcounter{equation}{0}

In this capter we shall demonstrate that non-interacting q-deformed systems
can be viewed as non-deformed systems with complicated, momentum dependent
interactions. In other words interacting systems can be described by a "free"
system based on a q-deformed kinematics.

We will start with the algebra SU$_q$(1,1) with its generators expressed
in terms of the generators of the undeformed algebra. Then we use the fact
that the undeformed generators of SU(1,1) can be represented in terms of the
creation and annihilation operators of the usual harmonic oscillator
(see (5.5)).\\[2mm]
This will give a manifest expression for the Hamiltonian with a complicated
momentum dependent interaction.\\

{}From the work of Curtright and Zachos \cite{curt} we learn how to express the
generators $T_+$, $T_-$, $T_3$ from the algebra (3.2) in terms
of the undeformed generators\footnote{$[x]_q = (q^x-q^{-x})/(q-q^{-1})$}:
\bea
T_+ &=& \sqrt{2}\,q^{-j_0}\,\sqrt{\frac{[j_0+j]_q [j_0-j-1]_q}{(j_0+j)
  (j_0-j-1)}}\;\;j_+ \nn\\
T_- &=& \sqrt{2}\,q^{-j_0}\,j_-\;\sqrt{\frac{[j_0+j]_q [j_0-j-1]_q}{(j_0+j)
  (j_0-j-1)}} \\
T_3 &=& \frac{1}{\lambda}\left( 1 - q^{-4 j_0} \right) \nn
\eea
where\\[-3mm]
\bea
j_+ j_- - j_- j_+ &=& j_0 \nn\\
j_0 j_{\pm} - j_{\pm} j_0 &=& \pm j_{\pm}
\eea
and $j$ is defined through the Casimir:
\be
{\cal C}_0 = j(j+1) = 2j_+j_- + j_0(j_0-1) = 2j_-j_+ + j_0(j_0+1) \nn
\ee
\\
For the compact version SU$_q$(2) these formulas are easy to obtain from the
representations of $T_+$, $T_-$ and $T_3$. As the algebraic relations are
the same for SU$_q$(1,1) (only the conjugation properties are different)
we can use these formulas to express $W_+$, $W_-$ and $W_0$ from (2.12):
\bea
W_+ &=& \sqrt{\frac{2r}{r+r^{-1}}}\; r^{-j_0}\; z_r^{\frac{1}{2}}\;j_+\nn\\
W_- &=& \sqrt{\frac{2r}{r+r^{-1}}}\;j_- \;z_r^{\frac{1}{2}}\; r^{-j_0}\\
W_0 &=& \frac{1}{r-r^{-1}}\; \big( 1-r^{-2j_0}\,\frac{r^{2j+1}+r^{-(2j+1)}}
  {r+r^{-1}} \big) \nn\\
z_r &=& \frac{[j_0+j]_r [j_0-j-1]_r}{(j_0+j)(j_0-j-1)}\nn
\eea
\\[-1mm]
The undeformed generators can now be represented through the creation and
annihilation operators of the undeformed harmonic oscillator:
\bea
j_0 &=& \frac{1}{4}\;(a^2-(a^+)^2) \nn\\
j_+ &=& \frac{1}{4\sqrt{2}} \; (a+a^+)^2 \\
j_- &=& -\frac{1}{4\sqrt{2}} \; (a-a^+)^2 \nn
\eea
In this representation we find $j(j+1) = -\frac{3}{16}$ , as expected. \\[2mm]
To make the transition from the W generators to the A,B,C generators we have
to find an expression for u. We use formulas (2.2) and (3.6) and obtain:
\be
u = q^{-\frac{1}{2}} q^{2 j_0}
\ee
Here we have made the identification $r=q^2$, as it is necessary for the
identification of the algebras (see (2.13)). It now leads to the following
expressions:
\bea
A &=& \alpha\, \qiq\, \sqrt{\frac{2q}{q^2+\qiq}}\;j_-\,z_{q^2}^{\frac{1}{2}}
\nn\\
B &=& \beta\, \sqrt{\frac{2q}{q^2+\qiq}}\;z_{q^2}^{\frac{1}{2}}\,j_+ \\
C &=& \gamma\,q^{-\frac{1}{2}}\,\frac{1}{q^2-q^{-2}}\:\left(q^{2j_0} -
q^{-2j_0}\,\frac{q+q^{-1}}{q^2+q^{-2}} \right) \nn \\
\alpha \beta &=& \frac{1}{q^3}(1+q^2)^2(1+q^4) \hspace{2cm}
\gamma = \frac{i}{q}\,(1+q^2)(1+q^4) \nn
\eea
$A$ and $B$ are hermitian operators. $\frac{1}{2} B$ is the free
Hamiltonian in the q-deformed
phase space, it was diagonalized in chapter 4.\\[2mm]
In a more physical
representation the Hamiltonian can be expressed in terms of the undeformed
phase space variables $x$ and $p_x = -i\frac{\partial}{\partial x}$,
a natural choice for $\beta$ is $ \beta\,\sqrt{2q/(q^2+q^{-2})} = 2 \sqrt{2}$:
\bea
{\cal H} \:=\: \frac{1}{2}\,B &=& 2\;\; \mbox{\large{\it p$_x$}} \;\;
\sqrt{\frac{q+q^{-1} - 2 \cos((x p_x+p_x x)\,h)}{(q^2-q^{-2})^2 \,
((x p_x+p_x x)^2+1)}}\;\;\: \mbox{\large{\it p$_x$}} \nn\\
q &=& e^h
\eea
The hermiticity of this Hamiltonian is explicit, because $q+q^{-1} \geq 2$ for
$q$ positive.\\
Observe that the operators are not really in the denominator (a physical way
of speaking).\\[1mm]
The first terms in an expansion in $h$ are
\be
{\cal H}\, = \,\frac{1}{2}\; p_x\;\left\{ 1- \frac{1}{2}\,h^2\,\big(
\frac{5}{4}+
\frac{2}{4!}(x p_x + p_x x)^2 \big) + \cdots \right\}\;\; p_x
\ee

With the previous normalization the generators $A$ and $B$ are
invariant under $q \rightarrow q^{-1}$ and thus are functions of $h^2$ only.
An explicit expansion of (6.7) and (6.6) to first order in $h$ leads to
our result (5.18) obtained by perturbation theory.
\bea
A &=& A_0 \,=\, 2\sqrt{2}\;j_- \nn\\
B &=& B_0 \,=\, 2\sqrt{2}\;j_+ \nn\\
C &=& C_0 + \frac{3}{2}(C_0+i)\,h \,=\, 4i\,j_0 + \frac{3}{2}(4i\,j_0+i)\,h\\
u &=& 1 + \frac{i}{2}(i-C_0)\,h\nn
\eea
As $u$ depends on $h$ linearly, the commutator of $u$ and $B$ will depend
linearly on $h$ as well. To first order in $h$ it is
\be
[\,u\, ,\,B\,] \,=\, 2\,h\,B
\ee
This generates the first order correction to the spectrum of $B$, it is
linear in $h$ and consistent with the exact result (4.3). This now justifies
our assumption $X_1=0$ and Im$<n-2|K_1^+|n>\,=\,0$ in the perturbative
treatment of the algebra.

Finally let us give the Hamiltonian for the q-deformed harmonic oscillator
in terms of the undeformed canonical variables:
\bea
{\cal H} &=& \frac{1}{2}\,(A + B) \nn\\
         &=& 2 \;\, \mbox{\large{\it p$_x$}}\;\;
         \sqrt{\frac{q+q^{-1} - 2 \cos((x p_x+p_x x)\,h)}
      {(q^2-q^{-2})^2 \,((x p_x+p_x x)^2+1)}}\:\;\;\mbox{\large{\it p$_x$}} \\
         & &   +\, \frac{1}{2q^2} (1+q^2)^2
        \;\;\mbox{\large{\it x}}\;\;
          \sqrt{\frac{q+q^{-1} - 2 \cos((x p_x+p_x x)\,h)}
          {(q^2-q^{-2})^2 \,((x p_x+p_x x)^2+1)}}\;\;\:\mbox{\large{\it x}} \nn
\eea
\\

\subsection*{Acknowledgement}

We would like to thank Prof. M. Marinov for stimulating discussions that
initiated this work.\\
We also would like to thank H. Ewen and A. Ruffing for helpful discussions.

\end{document}